\documentclass[12pt]{article}
\title{A complete diophantine characterization of the rational torsion of an elliptic curve}
\author{Irene Garc\'{\i}a--Selfa \and Jos\'e M. Tornero}
\date{March, 2007}

\newcommand{\QQ}{{\mathbf Q}}
\newcommand{\ZZ}{{\mathbf Z}}
\newcommand{\vs}{\vspace{.15cm}}
\newcommand{\VS}{\vspace{.3cm}}

\begin{document}

\maketitle

\abstract{We give a complete characterization for the rational torsion of an elliptic curve in terms of the (non--)existence of integral solutions of a system of diophantine equations.}

\vs\vs


\vs


\section{Introduction}

In this paper we consider elliptic curves defined over $\QQ$. As it is known \cite{Cassels, Silverman}, each one of such curves is birationally equivalent to one, say $E$; given by an equation of the type
$$
E: Y^2=X^3+AX+B, \quad \mbox{ with }A,B, \in \ZZ;
$$
called short Weierstrass form, where it must hold $\Delta = 4A^4+27B^2 \neq 0$. The set of rational points of its projective clausure, noted $E(\QQ)$, is a finitely generated abelian group (Mordell--Weil Theorem \cite{Mordell, Weil}) and its torsion part, noted $T(E(\QQ))$ has been exhaustively described by Mazur \cite{Mazur, MazurIHES} as isomorphic to one of the following groups:
$$
\begin{array}{lcl}
\ZZ / n \ZZ & \mbox{ for } & n=1,2,...,10,12; \\
\ZZ / 2 \ZZ \times \ZZ / 2n \ZZ & \mbox{ for } & n=1,2,3,4.
\end{array}
$$

We will give a procedure for finding out which of these groups is the actual torsion subgroup of our curve, by characterizing the existence of points of order $n$ in terms of the existence of an integral solution for a system of diophantine equations. 

Note that only orders which are prime or pure prime powers have to be considered; for the rest of the cases may be solved by joining systems, according to the factorization of the order we are interested on.

Although the actual system depends on $n$, we will see that they are fairly similar systems, in more than one sense. The main result can be stated as follows:

\vs

\noindent {\bf Theorem.--} Let $E:Y^2=X^3+AX+B$ be an elliptic curve, with $A,B \in \ZZ$. For every $n \in \{3,4,5,7,8,9\}$ there are, at most, $4$ quasi--homogeneous polynomials $P_n,Q_n,R_n,S_n \in \ZZ [z_1,...,z_4]$ such that $E$ has a rational point of order $n$ if and only if there exists an integral solution for the system
$$
\Sigma_n : \left\{ 
\begin{array}{rcl} 
P_n(z_1,...,z_4) &=& 6^2 \cdot A\\ 
Q_n(z_1,...,z_4) &=& 6^3 \cdot B \\
R_n(z_1,...,z_4) &=& 0 \\
S_n(z_1,...,z_4) &=& 0 \\
\end{array} \right.
$$

\vs

\noindent {\bf Remark.--} The case $n=2$ will be treated separatedly, as not only the existence, but the number of solutions is important for determining whether the torsion group is cyclic or not. However, we will get systems of the above type anyway, which we will call $\Sigma^{(1)}_2$ and $\Sigma^{(3)}_2$, for detecting the existence of one or three points of order two, respectively.

\noindent {\bf Remark.--} Actually, for most cases, less equations and variables have to be considered, according to the following table:

\begin{tabular}{llc}
\\
{\bf Case} & {\bf Equations and variables} & {\bf Max. degree} \\
\hline \\
$n=3$ & $2$ equations in $2$ variables & $6$ \\
$n=4$ & $2$ equations in $2$ variables & $3$ \\
$n=5$ & $3$ equations in $3$ variables & $4$ \\
$n=7$ & $4$ equations in $4$ variables & $4$ \\
$n=8$ & $3$ equations in $3$ variables & $6$ \\
$n=9$ & $3$ equations in $3$ variables & $10$\\
\\
\end{tabular}

The remarkable equality between the number of equations and variables may be of some help if we try to use some numerical algorithm to find integral solutions.

\noindent {\bf Remark.--} Now and then we will use the Nagell--Lutz Theorem \cite{Nagell, Lutz, Cassels, Silverman} which, among more specific results, states that the coordinates of torsion points must be integers, under our assumptions.

\vs

There exist some previous descriptions of torsion structures in terms of diophantine equations, notably by Ono \cite{Ono}, for the non--cyclic case and by Qiu and Zhang \cite{QZ,QZ2} for the even cyclic case (strongly influenced by Ono's ideas). Both papers took advantage of the possibility of choosing $(0,0)$ to be a $2$--torsion point, which is not suitable in the general case. Instead, while trying to cope with the odd order case we eventually arrived to a global strategy for finding equations which led to the previous theorem. 

Both Ono's and Qiu--Zhang's results have been used to study the behaviour of torsion structures when a bigger ground field is considered (typically a quadratic extension of $\QQ$), for instance in \cite{Kwon, Fujita, Fujita2}. We hope that our characterization will shed some light also in this matter, where much is yet to be known.

As a final remark, we should note that the resulting systems seem to be suitable for a computational diophantine attack, either classical and/or modern (see for instance \cite{DEQ, Smart} for a good account on each approach). However, it seems difficult to us that this strategy can beat the existing algorithms \cite{GOT} which can compute torsion subgroups of curves with discriminant $\Delta \sim 2^{20000}$ in only a few minutes.

\section{Points of order $3$ and $9$}

Assume our curve $E$ has a point $(x,y)$ of order three. As it is known \cite{Silverman}, $(x,y) \in E(\QQ)$ is a point of order three if and only if
$$ 
\Psi_3 (x) = 3x^4+6Ax^2+12Bx-A^2 = 0. 
$$

If $(x,y)=(0,z_2)$, then we must have $A=0$ and therefore, $B=z_2^2$.
So we will suppose $x\neq 0$ from now on, and divide $\Psi_3 (x)$ by $x^4$ to
obtain the conic 
$$
3+6X+12Y-X^2=0,
$$ 
where $X=A/x^2$ and $Y=B/x^3.$ Thus we can parametrize the conic by
$$ 
X = \frac{3(t+1)}{t-3},\ \, Y =\frac{-t^2+6t+3}{(t-3)^2},\ \   t \in \QQ \setminus \{3\}
$$

Now if we set $t=p/q$ for $p,q \in \ZZ$ with $\gcd(p,q)=1$ and $q\neq 0, \, p-3q \neq 0$,
then we have
$$  
\frac{A}{x^2}=\frac{3(p+q)}{p-3q},\ \  \frac{B}{x^3}=\frac{-p^2+6pq+3q^2}{(p-3q)^2}, 
$$

Writing $u=p-3q$ we get 
$$
A=3(u+4q)x^2/u,\ \ B=(12q^2-u^2)x^3/u^2.
$$ 

Thus,
$$ 
y^2=\frac{3(u+2q)^2}{u^2}x^3, 
$$
and, since $(x,y) \in E(\QQ)$, there exist $z_1,z_2 \in \ZZ$ such that
$x=3z_1^2$ and $y=\pm(9z_1^3+z_2)$ with $z_2=18z_1^3q/u$. 
Therefore,
$$ 
A=27z_1^4+6z_1z_2, \ \ B=z_2^2-27z_1^6,
$$
which also fits the case $x=0$.

Hence we have found:
$$
\Sigma_3: \left\{ \begin{array}{rcrcl}
6^2 \cdot A &=& P_3 &=& \displaystyle 6^2 \left( 27z_1^4+6z_1z_2 \right) \\ 
6^3 \cdot B &=& Q_3 &=& \displaystyle 6^3 \left( z_2^2-27z_1^6 \right)
\end{array} \right.
$$

\noindent {\bf Remark.--} The polynomials obtained are, in fact, quasi--homogeneous. If we get to compute the solutions for the systems, we obtain all points of order three, which are
$$
\{ (3z_1^2,9z_1^3+z_2), \; (3z_1^2,-(9z_1^3+z_2)) \}.
$$

Conversely, we can easily check that $(3z_1^2,\pm (9z_1^3+z_2))$ has
order three in the curve $ Y^2 = X^3 + (27z_1^4+6z_1z_2)X +(z_2^2-27z_1^6)$. This shows that the conditions are also sufficient. 

\vs\vs

Assume now that $E$ has a rational point of order nine. Then there must be rational points
of order three in the curve, thus we already have 
$$
\begin{array}{rcl}
P_9 &=& \displaystyle 6^2(27z_1^4+6z_1z_2) \\ 
Q_9 &=& \displaystyle 6^3(z_2^2-27z_1^6)
\end{array}
$$

Also, if $(x,y)$ is a point of order nine, the first coordinate of $[3](x,y)$ must be  $3z_1^2$, as $[3](x,y)$ has order three.

Using the chord--tangent and duplication formulae, the first coordinate in $[3](x,y)$ is
$$
3z_1^2= \frac{N(A,B,x)}{(3x^4+6Ax^2+12Bx-A^2)^2}.
$$
where
$$
N (A,B,x) = x^9-12Ax^7-96Bx^6+30A^2x^5-24ABx^4+(36A^3+48B^2)x^3
$$
$$
+48A^2Bx^2+(96AB^2+9A^4)x+8B(A^3+8B^2).
$$

So $(x,y) \in E(\QQ)$ has order nine if and only if
$$
R_9 = x^9-12Ax^7-96Bx^6+30A^2x^5-24ABx^4+(36A^3+48B^2)x^3+48A^2Bx^2
$$
$$
+(96AB^2+9A^4)x+8B(A^3+8B^2)-3z_1^2(3x^4+6Ax^2+12Bx-A^2)^2=0.
$$

For writing it our way, we should make the appropriate substitutions 
$$
A = 27z_1^4+6z_1z_2, \qquad B=z_2^2-27z_1^3; 
$$ 
to obtain a quasi--homogeneous polynomial, for instance with weights $w(z_1)=1, 
\, w(z_2) = 3, \, w(x)=2$; weighted degree $18$ (with the previous choice) and $69$ monomials. This polynomial is the most complicated we have to deal with, but it turns out it can be written in a more amusing way. If we perform the Tchirnhausen transformation w.r.t. $x$; that is, we make
$$
z_3 = x - 3z_1^2,
$$ 
we find out our polynomial $R_9$ can be written as
$$
\begin{array}{rcl} 
R_9 &=& z_3^9 - 72z_1(9z_1^3+z_2) z_3^7 - 48(27z_1^3+2z_2)(9z_1^3+z_2)z_3^6 \\
&& \quad - 864z_1^2(9z_1^3+z_2)^2 z_3^5 -144z_1(9z_1^3+z_2)^3z_3^4 \\
&& \quad + 48(45z_1^3+z_2)(9z_1^3+z_2)^3 z_3^3 + 1728 z_1^2(9z_1^3+z_2)^4z_3^2 \\
&& \quad +576z_1(9z_1^3+z_2)^5 z_3 + 64 (9z_1^3+z_2)^6
\end{array}
$$
which, taking into account the given parametrization of points of order three, proves the following funny fact:

\vs

\noindent {\bf Corollary.--} Let $(x,y)$ be a point of order nine, and write $[3](x,y) = (x',y')$. Then $(x-x')|(2y')^6$.

\vs

In order to get the easiest possible system we will perform
$$
z_2 \longmapsto z_2 - 9z_1^3
$$
and our system is then
$$
\Sigma_9: \left\{ \begin{array}{rcrcl}
6^2 \cdot A &=& P_9 &=& \displaystyle 6^2 \left( 27z_1^4+6z_1z_2 \right) \\ 
6^3 \cdot B &=& Q_9 &=& \displaystyle 6^3 \left( z_2^2-27z_1^3 \right) \\
0 &=& R_9 &=& \displaystyle z_3^9 - 72z_1z_2 z_3^7 - 48 \left( 9z_1^3+2z_2 \right )z_2 z_3^6 \\
&& && \displaystyle - 864z_1^2z_2^2 z_3^5 -144z_1z_2^3z_3^4 + 48 \left( 36z_1^3+z_2 \right)z_2^3 z_3^3 \\
&& && + 1728 z_1^2z_2^4z_3^2 +576z_1z_2^5 z_3 + 64 z_2^6
\end{array} \right.
$$

\noindent {\bf Remark.--} If one wants to use this system for deciding whether a given curve has a torsion group of order nine, it is wise to solve first 
$$
6^2 \cdot A = P_9, \qquad 6^3 \cdot B = Q_9,
$$
which is a system of two equations in two variables. One then moves on to solving $R_9$ as a polynomial in $z_3$. Notice then that this polynomial, which has degree nine, can only have three different integer roots corresponding to the first coordinates of points of order nine, since an elliptic curve can only have six such different points, and, in Weierstrass short form, three of them should be symmetric to the other three w.r.t. the $x$--axis.

\section{Points of order $5$ and $7$}

Orders five and seven share an ad--hoc strategy for obtaining the respective systems $\Sigma_5$ and $\Sigma_7$. Assume $(x_1,w_1) \in E(\QQ)$ has order
five. Then the first coordinate of $[4](x_1,w_1)$ must be again $x_1$, and if we write  $[2](x_1,w_1) = (x_2,w_2)$, then it must hold $x_1 \neq x_2$. We
will have, from duplication formula,
$$ 
\frac{x_1^4-2Ax_1^2-8Bx_1+A^2}{4(x_1^3+Ax_1+B)} = x_2;\ \ \frac{x_2^4-2Ax_2^2-8Bx_2+A^2}{4(x_2^3+Ax_2+B)}= x_1.
$$

This equations yield a system whose (naive) general solution for $B$ is given by
$$
\left\{ x_1=x_1,\ x_2=x_2,\ A=A,\ B=\frac{1}{4}(x_1+x_2)(x_1^2-4x_1x_2+x_2^2-2A) \right\},
$$
where $x_1,x_2,A \in \ZZ$ and $A$ must verify the second degree polynomial
$$
A^2+(2x^2+2xv+2v^2)A-x^4+x^3v+9x^2v^2+xv^3-v^4=0.
$$ 

The discriminant of this equation is 
$$
4(x_1-x_2)^2(2x_1+x_2)(x_1+2x_2),
$$
hence there must exist $x_3 \in \ZZ$ such that 
$$
x_3^2=(2x_1+x_2)(x_1+2x_2).
$$

With these notations we can finally write
$$ 
A = -x_1^2-x_1x_2-x_2^2 + (x_1-x_2)x_3,
$$
$$ 
B = \frac{-1}{4}(x_1+x_2)(-3x_1^2+2x_1x_2-3x_2^2+
2(x_1-x_2)x_3).
$$

Now, since 
$$
w_1^2=\frac{(x_1-x_2)^2(3x_1+2x_3+3x_2)}{4},
$$ 
there must be $x_4 \in \ZZ$ verifying 
$$
x_4^2=3x_1+2x_3+3x_2.
$$

This completes a system $\Sigma'_5$, given by
$$
\begin{array}{rcl}
A &=& \displaystyle \left(-x_1^2-x_1x_2-x_2^2 + (x_1-x_2)x_3 \right) \\
4B &=& -(x_1+x_2)(-3x_1^2+2x_1x_2-3x_2^2+2(x_1-x_2)x_3) \\
0 &=& x_3^2-(2x_1+x_2)(x_1+2x_2) \\
0 &=& x_4^2-(3x_1+2x_3+3x_2).
\end{array}
$$

Note that the polynomials obtained in this case are actually homogeneous, except for $S_5$. We can also get the torsion points with a little extra work. With our previous notations it must hold 
$$
w_2^2=(x_1-x_2)^2(3x_1-2x_3+3x_2)/4,
$$ 
hence we can find $t\in \ZZ$ such that with 
$$
t^2=3x_1-2x_3+3x_2.
$$ 

Then the full list of torsion points of order five is
$$
\left\{ \left(x_1, \, \pm \frac{(x_1-x_2)x_4}{2} \right), \; \left( x_2, \, \pm \frac{(x_1-x_2)t}{2} \right) \right\}.
$$

By plugging these points into the curve we see that conditions given by $\Sigma'_5$ are sufficient to guarantee the existence of torsion points of order five. We then eliminate $x_3$ from the last equation and rename
$$
z_1 = x_1, \quad z_2 = x_2, \quad z_3 = x_4
$$
to finally get the equivalent system $\Sigma_5$:
$$
\Sigma_5: \left\{ \begin{array}{rcrcl}
6^2 \cdot A &=& P_5 &=& \displaystyle 18 \cdot \displaystyle \left[ (z_1-z_2)^2-6z_1^2+(z_1-z_2)z_3^2 \right] \\
6^3 \cdot B &=& Q_5 &=& \displaystyle 54 \cdot \displaystyle \left( z_1+z_2 \right) \left[ 6z_1^2-2z_1z_2 +(z_2-z_1)z_3^2 \right] \\
0 &=& R_5 &=& \displaystyle z_3^4-6z_3^2 \left( z_1+z_2 \right) + \left( z_1-z_2 \right)^2 \\
\end{array} \right.
$$

\vs\vs

Let us do then the case $n=7$. Assume $(x_1,w_1) \in E(\QQ)$ is a point of order seven, and denote $[2](x_1,w_1)=(x_2,w_2)$ and $[3](x_1,w_1)=(x_3,w_3)$, being $x_1,$ $x_2$ and $x_3$ different integers. From the duplication formula, as above, we have
$$ 
F_1(x_1,x_2) = x_1^4-2Ax_1^2-8Bx_1+A^2 - 4(x_1^3+Ax_1+B)x_2 = 0. 
$$

From the expression of $x_3$ used in the order nine case we can write
$$ 
F_2(x_1,x_3) = x_1^9-12Ax_1^7-96Bx_1^6+30A^2x_1^5-24ABx_1^4+
$$
$$
+(36A^3+48B^2)x_1^3+48A^2Bx_1^2 +(96AB^2+9A^4)x_1+8B(A^3+8B^2)-
$$
$$
- x_3(3x_1^4+6Ax_1^2+12Bx_1-A^2)^2 = 0. 
$$

Also, using the addition formula for $(x_1,w_1)$ and $(x_2,w_2)$, we can
also get $x_3$ as an expression in $x_1$ and $x_2$. Now, the
difference between our two expressions for $x_3$ (one in $x_1$ and
another in $x_1$ and $x_2$) must vanish, and we get
$$
F_3(x_1,x_2)=0,
$$
where $F_3$ is a polynomial in $x_1$ and $x_2$ with $139$ terms whose degree is $20$.

We will get still another relationship between our parameters,
using now two different ways to express $x_3$, as it is also the first coordinate of $[4](x_1,w_1)$. The first one comes from $[2](x_2,w_2)$ and depends only on $x_2$, and the second one from the addition formula for $(x_1,w_1)$ and $(x_3,w_3)$ which
depends on $x_1$ and $x_3$. Thus we get 
$$
F_4(x_1,x_2,x_3)=0,
$$
where $F_4$ is a polynomial in $x_1,x_2$ and $x_3$ with $321$ terms whose
degree is $18$.

\vs

Now if we consider the ideal $\langle F_1,F_2,F_3,F_4 \rangle \subset \QQ[x_1,x_2,x_3]$ and calculate a Gr\"obner basis with respect to the graded reverse lexicographic order we get a basis consisting of five polynomials. The general solution for $\{A,B\}$ of this system is
$$
\Big\{ x_1=x_1,\ x_2=x_2,\ x_3=x_3,\ A=A,
$$
$$
B=\frac{1}{4} \left( x_3x_1^2-2x_1x_2x_3+x_2^2x_3+x_1^3-4x_2x_1^2-x_2^2x_1-2(x_2+x_1)A \right) 
\Big\},
$$
where $x_1,x_2,x_3,A \in \ZZ$ and $A$ verifies
$$
A^2+(2x_1^2+2x_2^2+2x_2x_1)A-
$$
$$
x_1^4+3x_2x_1^3-2x_1^3x_3+x_2^3x_1+3x_2x_3x_1^2+6x_1^2x_2^2-x_2^3x_3=0.
$$

As a side remark, it was not easy to find this relation, which parallels case $n=5$. In fact, for most methods (including most orders considered for Gr\"obner bases computations), $A$ ended up verifying a fourth degree equation, which became then a dead end for our hopes (see below).

The discriminant of the previous second degree equation verified by $A$ is
$$
4(x_2+2x_1)(x_1-x_2)^2(x_1+x_3+x_2)
$$ 
hence there must be $x_4 \in \ZZ$ with
$$
x_4^2=(x_2+2x_1)(x_1+x_3+x_2).
$$ 

Therefore
$$
A=-x_1^2-x_2^2-x_1x_2+ x_4(x_1-x_2), 
$$
$$
4B=(3x_1^3+x_3x_1^2+3x_2^2x_1+x_2^2x_3-2x_1x_2x_3+2x_2^3+ 2(x_2^2-x_1^2)x_4). 
$$

Now again, as $(x_1,w_1)$ has order 7, we can write $x_3$ as the first coordinate of $[3](x_1,w_1)$ and $[4](x_1,w_1)$, to get $A$ as a root of the polynomial
$$
Z^4+6x_2(x_1+x_2)Z^3+(-8x_1^4+15x_1^2x_2^2+11x_2x_1^3+23x_2^3x_1+13x_2^4)Z^2
$$
$$
+(-24x_2x_1^5+58x_2^3x_1^3+42x_2^4x_1^2-6x_2^2x_1^4+4x_1^6+24x_2^5x_1+10x_2^6)Z
$$
$$
-(x_1^2+x_1x_2+x_2^2)(x_2^6-3x_2^5x_1-33x_2^4x_1^2-20x_2^3x_1^3+39x_2^2x_1^4-12x_2x_1^5+x_1^6),
$$
which was the typical result we mentioned above.

Hence, if we substitute $Z$ by the above expression of $A$, we get
$$
-x_4^4 + 2(2x_1+x_2)x_4^3+(2x_1+x_2)(x_1-x_2)x_4^2 -
$$ 
$$
-2(2x_1+x_2)^3x_4+3(x_1^2+x_1x_2+x_2^2)(2x_1+x_2)^2=0.
$$

Since $x_4^2=(x_2+2x_1)(x_1+x_3+x_2)$, 
we have
$$
(2x_1+x_2)^2(3x_1^2-x_1x_3+x_1x_2-x_3^2-3x_2x_3+x_2^2-2(x_1-x_3)x_4) =0.
$$

It is an easy fact that $2x_1+x_2\neq 0$. Otherwise, $x_4=0$ and then it must hold
$$
x_1(5x_1+4x_3)-x_3(x_3-x_1) = 0,
$$
from the expression of $x_3$ above. We would then have 
$$
5x_1^2+5x_1x_3-x_3^2=0,
$$ 
which is impossible for different integers $x_1$ and $x_3$. Thus we obtain
$$
3x_1^2-x_1x_3+x_1x_2-x_3^2-3x_2x_3+x_2^2-2(x_1-x_3)x_4=0.
$$

Finally, since $(x_1,w_1) \in E(\QQ)$ and 
$$
w_1^2=\frac{1}{4}(x_1-x_2)^2(3x_1+2x_2+2x_4+x_3),
$$ 
there must be $x_5\in \ZZ$ such that 
$$
x_5^2=(3x_1+2x_2+2x_4+x_3).
$$ 

Thus we have proved that the necessary conditions for having a point of order seven lead yield the existence of solution for $\Sigma'_7$, with:
$$
\begin{array}{lcl}
A &=& (-x_1^2-x_2^2-x_1x_2+ x_4(x_1-x_2)) \\
4B &=& (3x_1^3+x_3x_1^2+3x_2^2x_1+x_2^2x_3-2x_1x_2x_3+2x_2^3+ 2(x_2^2-x_1^2)x_4) \\ 
0 &=& x_4^2-(x_2+2x_1)(x_1+x_3+x_2) \\
0 &=& 3x_1^2-x_1x_3+x_1x_2-x_3^2-3x_2x_3+x_2^2-2(x_1-x_3)x_4 \\
0 &=& x_5^2-(3x_1+2x_2+2x_4+x_3)
\end{array}
$$

Since also $w_2,w_3 \in \ZZ$, there will be $s,t \in \ZZ$ verifying 
$$
\begin{array}{rcl}
s^2 &=& 3x_1+2x_2+x_3-2s \\ 
t^2 &=& 4x_3^3-3x_3x_1^2-3x_2^2x_3-6x_1x_2x_3+4x_3x_1s-4x_3x_2s \\
&& \quad +3x_1^3+3x_1x_2^2-2sx_1^2+2x_2^3+2sx_2^2
\end{array}
$$

Now one can compute the coordinates of all order seven points and the curve is easily checked to have the desired torsion group. In order to arrive at $\Sigma_7$ we make
$$
z_1=x_1, \quad z_2 = x_2, \quad z_3 = x_5, \quad z_4 = x_4,
$$
and eliminate $x_3$ using the last equation from $\Sigma'_7$, to obtain
$$
\Sigma_7: \left\{ \begin{array}{rcl}
6^2 \cdot A &=& \displaystyle -6^2 \cdot \left[ (z_1+z_2)^2-z_1z_2 + z_4(z_2-z_1) \right] \\
6^3 \cdot B &=& \displaystyle 54 \cdot \left[ z_3^2 (z_1-z_2)^2 +4z_2(z_1^2+z_1z_2+z_1z_4-z_2z_4) \right] \\
0 &=& \displaystyle \left( 2z_1+z_2+z_4 \right)^2 - z_3^2 \left( 2z_1+z_2 \right) \\
0 &=& \displaystyle -z_3^4 + \left( 5z_1+z_2+6z_4 \right)z_3^2 + \\
&& \displaystyle \quad \left( -3z_1^2+3z_2^2-8z_4^2-18z_1z_4 -6z_2z_4 \right) 
\end{array} \right.
$$

\noindent {\bf Remark.--} We wanted all of our systems to have the previous form. If one is not so concerned about that one can find a simpler (in principle) system. Clearly, from the third equation
$$
z_3^2 = \frac{(2z_1+z_2+z_4)^2}{2z_1+z_2},
$$
which we can substitute in the second and fourth equation. Note that, in the fourth equation, we can clear the resulting denominator, as
$$
(2z_1+z_2)^2 = (2x_1+x_2)^2 \neq 0.
$$

One should be concerned about $z_3$ being an integer. This might not be guaranteed if we perform the above substitution. However, as
$$
2z_1+z_2 = 2x_1+x_2 = \left( \frac{A+3x_1^2}{2w_1} \right)^2,
$$
using the duplication formula, the above substitution is valid. 

Finally, 
$$
\left\{ \begin{array}{rcl}
0 &=& A \displaystyle + \left[ (z_1+z_2)^2-z_1z_2 + z_4(z_2-z_1) \right] \\
0 &=& 4 B (2z_1+z_2) - \displaystyle \left[ (2z_1+z_2+ z_4)^2 (z_1-z_2)^2 \right. \\
&& \displaystyle \quad \left. +4z_2(2z_1+z_2)(z_1^2+z_1z_2+z_1z_4-z_2z_4) \right] \\
0 &=& \displaystyle - z_4^4 + 2 \left( 2z_1+z_2 \right)z_4^3 + (2z_1+z_2)(z_1-z_2) z_4^2 \\
&& \displaystyle \quad -2(2z_2+z_2)^3z_4 + 3\left( z_1^2+z_1z_2+ z_2^2\right)(2z_1+z_2)^2 
\end{array} \right.
$$

\section{Points of even order}

We will finish with a quick discussion on points of even order; which is the most treated case in the existing literature. We noted before that order two has a special feature: it is important to know exactly the number of points, not only to make sure one such point exists. However, we already have the following
$$
(\alpha,\beta) \in E \mbox{ has order two } \; \Longrightarrow \; \beta = 0 \mbox{ and } \alpha^3 + A \alpha + B = 0.
$$

Hence the existence (and the number) of points of order two is that of the rational roots of $X^3+AX+B$. If we call the roots $\{ \alpha, \beta, \gamma \}$ (whether they are rational or not), using Cardano formulae it is straightforward to show that
$$
B = \beta^2\gamma + \gamma \beta^2 = \beta \gamma (\beta+\gamma), \qquad A = -\beta^2-\gamma^2- \beta \gamma = \beta \gamma - (\beta + \gamma)^2.
$$

Now, $\alpha \in \QQ$ if and only if so does $\beta \gamma$ and $\beta +\gamma$. Hence the existence of a rational root is equivalent to the existence of an integral solution for:
$$
\Sigma_2^{(1)}: \left\{ \begin{array}{rcrcl} 
6^2 \cdot A &=& P^{(1)}_2 &=& \displaystyle 6^2 \cdot \left( z_1 -z_2^2 \right) \\
6^3 \cdot B &=& Q_2^{(1)} &=& \displaystyle 6^3 \cdot z_1z_2.
\end{array} \right.
$$

In this context, the points of order two are those given by
$$
-z_2, \frac{1}{2} \left(-z_2 \pm \sqrt{z_2^2-4z_1} \right);
$$
so the existence of three points of order two is characterized by the existence of an integral solution for
$$
\begin{array}{rcl}
A &=& z_1-z_2^2 \\
B &=& z_1z_2 \\
0 &=& z_3^2-z_2^2+4z_1
\end{array}
$$

Eliminating $z_1$ and renaming gives then:
$$
\Sigma_2^{(3)}: \left\{ \begin{array}{rcrcl}
6^2 \cdot A &=& P_2^{(3)} &=& \displaystyle -9 \cdot \left( z_1^2+3z_2^2 \right) \\
6^3 \cdot B &=& Q_2^{(3)} &=& \displaystyle 54 \cdot \left( z_2^3-z_2z_3^2 \right) \\
\end{array} \right.
$$

\VS

Moving on to order four, assume the point $(x_1,w_1) \in E(\QQ)$ has order four, what happens if and only if $[3](x_1,w_1) = (x_1,-w_1)$, and $w_1 \neq 0.$ Using the same expression as in order nine for the first coordinate of $[3]P$, we quickly 
get 
$$
B =\frac{1}{4} \left( 5x_1^3-Ax_1 \pm \sqrt{(3x_1^2-2A)(A+3x_1^2)^2} \right).
$$ 

Therefore, there must be $x_2 \in \QQ$ such that $3x_1^2-2A=x_2^2.$ Then we
have 
$$
A= \frac{1}{2}(3x_1^2-x_2^2),
$$ 
and 
$$
B=\frac{1}{8} (x_2+x_1)(7x_1^2+2x_1x_2-x_2^2).
$$

The second coordinate verifies
$$
w_1^2=\frac{1}{8}(3x_1-x_2)(x_2+3x_1)^2,
$$ 
so it must be $x_2\neq \pm 3x_1$ and there must be $x_3 \in \ZZ \setminus \{0\}$
such that 
$$
\frac{3x_1-x_2}{2}=x_3^2.
$$ 

Our provisional system $\Sigma'_4$ then will be
$$
\begin{array}{rcl}
2A &=& 3x_1^2-x_2^2 \\
8B &=& (x_2+x_1)(7x_1^2+2x_1x_2-x_2^2) \\
0 &=& 2x_3^2-3x_1+x_2
\end{array}
$$

Note that $x_2,x_3\in \ZZ,$ since $x_1,A \in \ZZ$, and $x_1$ and $x_2$ have the same parity; because $2A=3x_1^2-x_2^2$.

In order to prove that the existence of solutions for $\Sigma'_4$ suffices, we can compute the full list of points of order four, which are
$$
\left\{ \left(x_1,\pm \frac{x_3(3x_1+x_2)}{2} \right) \right\}
$$ 
and check are they are, in fact, points of order 4 in 
$$
E: Y^2 = X^3 + \frac{3x_1^2-x_2^2}{2} X + \frac{(x_2+x_1)(7x_1^2+2x_1x_2-x_2^2)}{8}.
$$ 

As for orders five and seven, we can now eliminate $x_2$ and rename
$$
z_1=x_1, \quad z_2 = x_3
$$
to get 
$$
\Sigma_4: \left\{\begin{array}{rcrcl}
6^2 \cdot A &=& P_4 &=& \displaystyle 6^2 \left( -3 z_1^2  + 6 z_1 z_2^2  - 2 z_2^4 \right) \\
6^3 \cdot B &=& Q_4 &=& \displaystyle 6^3 \left( 2z_1-z_2^2 \right) \left( z_1^2+2z_1z_3^2 - z_3^2 \right) \\
\end{array} \right.
$$

\VS

Let us finally do the case of order eight. Assume then that $E$ has a rational point $(v,t)$ of order eight. In this case, there must be rational points of order four in the curve, therefore we already have 
$$
\begin{array}{rcl}
6^2 \cdot A &=& \displaystyle 6^2 \left( -3 z_1^2  + 6 z_1 z_2^2  - 2 z_2^4 \right) \\
6^3 \cdot B &=& \displaystyle 6^3 \left( 2z_1-z_2^2 \right) \left( z_1^2+2z_1z_3^2 - z_3^2 \right) \\
\end{array}
$$ 
with $z_1,z_2 \in \ZZ$ such that $z_2\neq \pm 3z_1$, and $[2](v,t)=(z_1,w_1)$.

Now using the duplication formula, we
have the equation 
$$
v^4-4z_1v^3+(6z_1^2-12z_1z_2+4z_2^4)v^2+(-48z_1^2z_2^2+40z_1z_2^4-8z_2^6-4z_1^3)v 
$$
$$
z_1^4-48z_1^3z_2^2+64z_1^2z_2^4-28z_1z_2^6+4z_2^8=0,
$$
whose solutions are
$$
v=z_1+z_3(z_2+z_4)
$$
for some $z_3,z_4 \in \ZZ$ such that 
$$
z_3^2+z_4^2-3z_1 = 0, \qquad z_4^2-z_2(2z_3+z_2)=0.
$$

Furthermore, we get
$$
t=\pm z_3z_4(z_2+z_3+z_4).
$$

So an elliptic curve $E:\ Y^2=X^3+AX+B,$ with $A, B \in \ZZ$,
has a point $(v,t)$ of order 8 in $E(\QQ)$ if and only if the system $\Sigma'_8$, given by
$$
\begin{array}{rcl}
6^2 \cdot A &=& \displaystyle 6^2 \left( -3 z_1^2  + 6 z_1 z_2^2  - 2 z_2^4 \right) \\
6^3 \cdot B &=& \displaystyle 6^3 \left( 2z_1-z_2^2 \right) \left( z_1^2+2z_1z_3^2 - z_3^2 \right) \\
0 &=& z_3^2+z_4^2-3z_1 \\
0 &=& z_4^2-z_2(2z_3+z_2)
\end{array}
$$
has a solution $(z_1,...,z_4) \in \ZZ^4$. As previously, by computing all points of order eight, the condition is easily shown to be sufficient. The customary elimination ($z_1$ in this case) and renaming ($z_i \longmapsto z_{i-1}$) gives finally:
$$
\Sigma_8 : \left\{ \begin{array}{rcrcl}
6^2 \cdot A &=& P_8 &=& \displaystyle -12 \cdot \left( z_1^4 - 4 z_1^2z_2^2 + z_2^4 \right)\\
6^3 \cdot B &=& Q_8 &=& \displaystyle 8 \cdot \left(z_1^2-2z_2^2 \right) \left(2z_1^4 - 8z_1^2z_2^2 - z_2^4 \right) \\
0 &=& R_8 &=& \displaystyle z_3^2-z_1 \left( 2z_2+z_1 \right)
\end{array} \right.
$$

As another example of unexpected properties which can be deduced from this equations, when we eliminate $z_1$ we can prove the following result.

\vs

\noindent {\bf Corollary.--} Let $(x,y)$ be a point of order four. If there are points of order eight, then $3x$ can be written as sum of two squares.

\section{Appendix: Irredundancy of $R_n,S_n$}

We will produce a few examples to show why our polynomials are in fact necessary in order to assure the desired torsion structure. Of course this does not mean one cannot find simpler systems characterizing $T(E(\QQ))$. Note that orders two, three and four do not need such examples to be found. Similarly, it is clear that $R_9$ cannot be dropped, as it would mean that every curve with points of order three would have points of order nine.

\VS

\noindent \fbox{\bf Order five}

\vs 

\noindent {\bf Example.--} Consider 
$$
E:\, Y^2 = X^3 -92X +480.
$$

If we set $\Sigma_5$, we find that 
$$
z_1 = 2, \; z_2 = 4, \; z_3 = 6
$$
is a solution for the first two equations. However
$$
R_5 (2,4,6) = 4.
$$

Accordingly, $T(E(\QQ))$ is trivial.

\VS

\noindent \fbox{\bf Order seven}

\vs 

\noindent {\bf Example.--} Consider
$$
E: \, Y^2 = X^3 +X + 23.
$$

Then 
$$
z_1=1, \; z_2 = 2, \; z_3 = -2, \; z_4 = -8;
$$
is a solution for the first three equations of $\Sigma_7$, but
$$
S_7 (1,2,-2,-8) = -443.
$$

Again, the torsion subgroup of $E$ is trivial.

\vs

\noindent {\bf Example.--} The curve
$$
E: \, Y^2 = X^3-4X
$$
has $T(E(\QQ))= (\ZZ /2 \ZZ) \times (\ZZ / 2 \ZZ)$. If we only onsider the first, second and fourth equation from $\Sigma_7$ we find the solution:
$$
z_1=2, \; z_2=-2, \; z_3=z_4=0;
$$
but
$$
R_7 (2,-2,0,0) = 4.
$$

\VS

\noindent \fbox{\bf Order eight}

\vs 

\noindent {\bf Example.--} Let us take
$$
E: \ Y^2 = X^3-92X+480.
$$

The first two equations from $\Sigma_8$ admit the common solution
$$
z_1=3, \; z_2=6.
$$

However that leads to
$$
R_8 (3,6,z_3) = z_3^2-54,
$$
which has no integral solutions. As expected, $T(E(\QQ)) = \ZZ / 4 \ZZ$.

\end{document}